\documentclass[12pt]{article}

\usepackage[english]{babel}
\usepackage{amssymb,a4wide}
\usepackage{amsmath,xypic,epsf,mathrsfs,theorem}

\usepackage[all]{xy}
 
\newtheorem{theorem}{Theorem}[subsection]
\newtheorem{lemma}[theorem]{Lemma}

\newtheorem{corollary}[theorem]{Corollary}

\newtheorem{proposition}[theorem]{Proposition}

{\theorembodyfont{\upshape}
\newtheorem{definition}[theorem]{Definition}
\newtheorem{remark}[theorem]{Remark}
\newtheorem{example}[theorem]{Example}
\newtheorem{question}[theorem]{Question}}

\numberwithin{equation}{section}
\numberwithin{theorem}{section}

\newcommand{\diag}{\mathrm{diag}}
\newcommand{\Id}{\mathrm{id}}
\newcommand{\GL}{\mathrm{GL}}

\newcommand{\setmin}{\! \setminus \!}

\newcommand{\ep}{{\varepsilon}}
 
\newcommand{\cA}{{\mathcal A}}

\newcommand{\cU}{{\cal U}}

\newcommand{\cB}{{\cal B}}

\newcommand{\cO}{{\cal O}}
\newcommand{\cT}{{\cal T}}

\newcommand{\C}{{\mathbb C}}

\newcommand{\T}{{\mathbb T}}

\newcommand{\Cs}{{$C^*$-al\-ge\-bra}}

\newcommand{\sh}{{$^*$-ho\-mo\-mor\-phism}}

\newenvironment{proof}[1][Proof:]
{\begin{trivlist}\item[]\textbf{#1} }
{\hbox{}\nobreak\hfill\quad\hbox{$\square$}\end{trivlist}}

\begin{document}
 
\title{Properly infinite $C(X)$-algebras and $K_1$-injectivity}

\author{Etienne Blanchard, Randi Rohde and Mikael R\o rdam}
\date{}
\maketitle


\vspace{1cm}

\begin{abstract} \noindent We investigate if a unital $C(X)$-algebra is
  properly infinite when all its fibres are properly infinite. We show
  that this question can be rephrased in several different ways,
  including the question if every unital properly infinite \Cs{} is
  $K_1$-injective. We provide partial answers to these questions, and
  we show that the general question on proper infiniteness of 
  $C(X)$-algebras can be reduced to establishing proper infiniteness
  of a specific $C([0,1])$-algebra with properly infinite fibres.
  \end{abstract}

\section{Introduction}

\noindent The problem that we mainly are concerned with in this paper
is if any unital $C(X)$-algebra with properly infinite fibres is
itself properly infinite (see Section~\ref{sec:C(X)} for a brief
introduction to $C(X)$-algebras). An analogous study was carried out in the
recent paper \cite{hrw} where it was decided when $C(X)$-algebras, whose
fibres are either stable or absorb tensorially a given strongly self-absorbing
\Cs, itself has the same property. This was answered in the
affirmative in \cite{hrw} under the crucial 
assumption that the dimension of the space $X$ is
finite, and counterexamples were given in the infinite dimensional
case.

Along similar lines, Dadarlat, \cite{md}, recently proved that
$C(X)$-algebras, whose fibres are Cuntz algebras, are trivial under some
$K$-theoretical conditions provided that the space $X$ is finite
dimensional. 

The property of being properly infinite turns out to behave very
differently than the property of being stable or of absorbing a
strongly self-absorbing \Cs. It is relative easy to see
(Lemma~\ref{neighborhood}) that if a fibre $A_x$ of a $C(X)$-algebra $A$ is
properly infinite, then $A_F$ is properly infinite for some closed
neighborhood $F$ of $x$.  The
(possible) obstruction to proper infiniteness of the $C(X)$-algebra is
hence not local. Such an obstruction is also
not related to the possible complicated
structure of the space $X$, as we can show that a counterexample, if
it exists, can be taken to be a (specific) $C([0,1])$-algebra
(Example~\ref{ex1} and Theorem~\ref{equiv conditions}). 
The problem appears to be related with some rather
subtle internal structure properties of properly infinite \Cs s.

Cuntz studied purely infinite---and in the process also properly
infinite---\Cs s, \cite{Cu}, where he among many other things (he was
primarily interested in calculating the $K$-theory of his algebras
$\cO_n$) showed that any unital properly infinite \Cs{} $A$ is
$K_1$-surjective, i.e., the mapping $\cU(A) \to K_1(A)$ is onto; and
that any purely infinite simple \Cs{} $A$ is
$K_1$-injective, i.e., the mapping $\cU(A)/\cU^0(A) \to K_1(A)$ is
injective (and hence an isomorphism). He did not address the question
if any properly infinite \Cs{} is $K_1$-injective. That question has
not been raised formally to our knowledge---we do so here---but it does appear
implicitly, eg.\ in \cite{Ror:cuntz} and in \cite{tw}, 
where $K_1$-injectivity of properly infinite \Cs s has to be assumed.

Proper infiniteness of \Cs s has relevance for existence (or rather
non-existence) of traces and quasitraces. Indeed, a unital \Cs{}
admits a 2-quasitrace if and only if no matrix algebra over the \Cs{}
is properly infinite, and  a unital \emph{exact} \Cs{}
admits tracial state again if and only if no matrix algebra over the \Cs{}
is properly infinite.

In this paper we show that every properly infinite \Cs{} is
$K_1$-injective if and only if every $C(X)$-algebra with properly
infinite fibres itself is properly infinite. We also show that a
\emph{matrix algebra} over any such $C(X)$-algebra is properly
infinite. Examples of unital \Cs s $A$, where $M_n(A)$ is properly
infinite for some natural number $n \ge 2$ but where $M_{n-1}(A)$ is
not properly infinite, are known, see \cite{Ror:sums} and \cite{ro2}, but still
quite exotic. 

We relate the question if a given properly infinite \Cs{} is
$K_1$-injective to questions regarding homotopy of projections
(Proposition~\ref{prop:equiv}). In 
particular we show that our main questions are equivalent to the
following question: is any non-trivial projection in the first copy of
$\cO_\infty$ in the full
unital universal free product $\cO_\infty \ast \cO_\infty$ homotopic to any
(non-trivial) projection in the second copy of $\cO_\infty$? The
specific $C([0,1])$-algebra, mentioned above, is perhaps not
surprisingly a sub-algebra of
$C([0,1], \cO_\infty \ast \cO_\infty)$.

Using ideas implicit in Rieffel's paper, \cite{Rieffel}, we construct in
Section~\ref{sec:examples} a $C(\T)$-algebra $\cB$ for each \Cs{} $A$ and
for each unitary $u \in A$ 
for which $\diag(u,1)$ is homotopic to $1_{M_2(A)}$; and $\cB$
is non-trivial if $u$ is not homotopic to $1_A$. In this way we relate
our question about proper infiniteness of $C(X)$-algebras to a
question about $K_1$-injectivity. 

The last mentioned author thanks Bruce Blackadar for
many inspiring conversations on topics related to this paper.

\section{$C(X)$-algebras with properly 
infinite fibres} \label{sec:C(X)}
A powerful tool in the classification of $C^\ast$-algebras is the
study of their projections. A projection in a $C^\ast$-algebra is
said to be \emph{infinite} if it is equivalent to a
proper subprojection of itself, and it is said to be \emph{properly
  infinite} if it is equivalent to two mutually
orthogonal subprojections of itself. 

A projection which is not infinite is said to be \emph{finite}. A
unital $C^\ast$-algebra is said to be finite, infinite, or properly
infinite if its unit is finite, infinite, or properly infinite,
respectively. If $A$ is a \Cs{} for which 
$M_n(A)$ is finite for all positive integers $n$, then $A$ is \emph{stably
finite}.

In this section we will study stability properties of proper
infiniteness under (upper-semi-)continuous deformations using 
the \emph{Cuntz-Toeplitz algebra} which is defined as follows.
For all integers $n \geq
2$ the Cuntz-Toeplitz algebra $\mathcal{T}_n$ 
is the universal $C^\ast$-algebra generated by $n$ isometries $s_1,\dots,
s_n$ satisfying the relation $$s_1s_1^\ast + \dots + s_ns_n^\ast \leq
1.$$ 

\begin{remark} \label{rem:propinf}
A unital $C^\ast$-algebra $A$ is properly infinite if and only if
$\mathcal{T}_n$ embeds unitally into $A$ for some $n\geq 2$, in which
case $\cT_n$ embeds unitally into $A$ for all $n \ge 2$. 
\end{remark}

\noindent
In order to study deformations of such algebras, let us recall a few
notions from the theory of $C(X)$-algebras.

Let $X$ be a compact Hausdorff space and $C(X)$ be the $C^\ast$-algebra
of continuous functions on $X$ with values in the complex field $\mathbb{C}$. 

\begin{definition} 
A $C(X)$-algebra is a $C^\ast$-algebra $A$ endowed 
with a unital $^\ast$-homo\-morphism from $C(X)$ 
to the center of the multiplier $C^\ast$-algebra $\mathcal{M}(A)$ of~$A$. 
\end{definition}

\noindent
If $A$ is as above and $Y\subseteq X$ is a closed subset, then we put
$I_Y = C_0(X\setminus Y)A$, which is a closed two-sided ideal in
$A$. We set $A_Y = A/I_Y$ and denote the quotient map by $\pi_Y$.

For an element $a \in A$ we put $a_Y = \pi_Y(a)$, and if $Y$ consists
of a single point $x$, we will write $A_x, I_x, \pi_x$ and $a_x$ in the
place of $A_{\{x\}}, I_{\{x\}}, \pi_{\{x\}}$ and $a_{\{x\}}$,
respectively. We say that $A_x$ is the \emph{fibre} of $A$ at
$x$.

The function 
$$
x\mapsto\| a_x\| =\inf\{\|\, [1-f+f(x)]a\| : f\in C(X)\}
$$
is upper semi-continuous for all $a \in A$ (as one can see using the
right-hand side identity above). A $C(X)$-algebra $A$ is said to be
\textit{continuous}  
(or to be a \textit{continuous $C^\ast$-bundle over $X$}) if 
the function $x\mapsto\| a_x\|$ is actually continuous 
for all element $a$ in $A$. 

For any unital \Cs{} $A$ we let $\cU(A)$ denote the group of unitary
elements in $A$, $\cU^0(A)$ denotes its connected component containing
the unit of $A$, and $\cU_n(A)$ and $\cU_n^0(A)$ are equal to
$\cU(M_n(A))$ and $\cU^0(M_n(A))$, respectively. 

An element in a \Cs{} $A$ is said to be \emph{full} if it is not
contained in any proper closed two-sided ideal in $A$. 

It is well-known (see for example \cite[Exercise~4.9]{LLR}) that if
$p$ is a properly infinite, full projection in a \Cs{} $A$, then $e
\precsim p$, i.e., $e$ is equivalent to a subprojection of $p$, 
for every projection $e \in A$. 

We state below more formally three more or less well-known 
results that will be used frequently throughout
this paper, the first of which is due to Cuntz, \cite{Cu}.

\begin{proposition}[Cuntz] \label{prop:propinf1}
Let $A$ be a \Cs{} which contains at least one properly infinite, full
projection. 
\begin{enumerate}
\item Let $p$ and $q$ be properly infinite, full projections in
  $A$. Then $[p] = [q]$ in $K_0(A)$ if and only if $p \sim q$.
\item For each element $g \in K_0(A)$ there is a properly infinite,
  full projection $p \in A$ such that $g = [p]$.
\end{enumerate}
\end{proposition}

\noindent The second statement is a variation of the Whitehead lemma.

\begin{lemma}\label{lm:unitary}
Let $A$ be a unital $C^\ast$-algebra.
\begin{enumerate}
\item Let $v$ be a partial isometry in
$A$ such that $1-vv^*$ and $1-v^*v$ are properly infinite and full
projections.  
Then there is a unitary element $u$ in $A$ such that $[u]=0$ in
$K_1(A)$ and $v = uv^\ast v$, i.e., $u$ extends $v$. 
\item Let $u$ be a unitary element $A$ such that $[u]=0$ in
  $K_1(A)$. Suppose there exists a 
  projection $p \in A$ such that $\|up-pu\| < 1$ and $p$ and $1-p$ are
  properly infinite and full. Then $u$ belongs to $\mathcal{U}^0(A)$. 
\end{enumerate}
\end{lemma}

\begin{proof} (i). 
It follows from Proposition~\ref{prop:propinf1}~(i) that $1-v^\ast v \sim
1-vv^\ast$, so there is a partial isometry $w$ such that $1-v^\ast v =
w^\ast w$ and $1-v v^\ast = ww^\ast$. Now, $z=v + w$ is a unitary
element in $A$ with $zv^*v = v$. The projection $1-v^*v$ is properly
infinite and 
full, so $1 \precsim 1-v^*v$, which implies that there is an isometry
$s$ in $A$ 
with $ss^* \le 1-v^*v$. As $-[z]=[z^*]=[sz^*s^*+(1-ss^*)]$ in $K_1(A)$ (see
eg.\ \cite[Exercise 8.9 (i)]{LLR}), we see that $u =
z(sz^*s^*+(1-ss^*))$ is as desired.

(ii). Put $x= pup+(1-p)u(1-p)$ and note that $\|u-x\|<1$. It follows
that $x$ is invertible in $A$ and that $u \sim_h x$ in $\GL(A)$. Let
$x=v|x|$  be the polar decomposition of $x$, where $|x| =
(x^*x)^{1/2}$ and $v = x|x|^{-1}$ is unitary. Then $u \sim_h v$ in $\cU(A)$
(see eg.\ \cite[Proposition 2.1.8]{LLR}), and $pv=vp$. We proceed to
show that $v$ belongs to $\cU^0(A)$ (which will entail that $u$
belongs to $\cU^0(A)$). 

Write $v=v_1v_2$, where
$$v_1 = pvp+(1-p), \qquad v_2 = p + (1-p)v(1-p).$$
As $1-p \precsim p$ we can find a symmetry $t$ in $A$ such that
$t(1-p)t \le p$. As $t$ belongs to $\cU^0(A)$ (being a symmetry), we
conclude that $v_2 \sim_h tv_2t$, and one checks that 
$tv_2t$ is of the form $w +
(1-p)$ for some unitary $w$ in $pAp$. It follows that $v$ is homotopic
to a unitary of the form $v_0 + (1-p)$, where $v_0$ is a unitary in
$pAp$. We can now apply eg.\ \cite[Exercise 8.11]{LLR} to conclude
that $v \sim_h 1$ in $\cU(A)$. 
\end{proof}

\noindent We remind the reader that if $p,q$ are projections in a
unital \Cs{} $A$, then $p$ and $q$ are homotopic, in symbols $p \sim_h
q$, (meaning that they can be connected by a continuous path of
projections in $A$) if and only if $q = upu^*$ for some $u \in
\cU^0(A)$, eg.\ cf.\ \cite[Proposition 2.2.6]{LLR}. 

\begin{proposition} \label{prop:propinf2}
Let $A$ be a unital \Cs. Let $p$ and $q$ be two properly infinite, full
projections in $A$ such that $p \sim q$. Suppose that there exists a properly
  infinite, full projection $r \in A$ such that $p \perp r$ and $q
  \perp r$. Then $p \sim_h q$.
\end{proposition}

\begin{proof} Take a partial isometry $v_0 \in A$ such that $v_0^*v_0=p$ and
$v_0v_0^* = q$. Take a subprojection $r_0$ of $r$ such that $r_0$ and
$r-r_0$ both are properly infinite and full. Put $v= v_0+r_0$. Then
$vpv^* = q$ and $vr_0 = r_0 = r_0v$. Note that $1-v^*v$ and $1-vv^*$
are properly infinite and full (because they dominate the properly
infinite, full projection $r-r_0$). Use Lemma~\ref{lm:unitary}~(i) to
extend $v$ to a unitary $u \in A$ with $[u]=0$ in $K_1(A)$. Now,
$upu^*=q$ and $ur_0 = vr_0 = r_0 = r_0v=r_0u$. Hence $u \in \cU^0(A)$ by
Lemma~\ref{lm:unitary}~(ii), and so $p \sim_h q$ as desired.
\end{proof}

\begin{definition} \label{def:K1inj}
A unital \Cs{} $A$ is said to be \emph{$K_1$-injective}
if the natural mapping
$$\cU(A)/\cU^0(A) \to K_1(A)$$
is injective. In other words, if $A$ is $K_1$-injective, and if $u$ is
a unitary element in $A$, then $u \sim_h 1$ in $\cU(A)$ if (and only
if) $[u]=0$ in $K_1(A)$. 
\end{definition}

\noindent One could argue that $K_1$-injectivity should entail that
the natural mappings $\cU_n(A)/\cU_n^0(A) \to K_1(A)$ be injective for
every natural number $n$. However there seem to be an agreement for
defining $K_1$-injectivity as above. As we shall see later, in
Proposition~\ref{prop:K1inj}, if $A$ is properly infinite, then the
two definitions agree. 

\begin{proposition}  \label{prop:pull-back}
Let $A$ be a unital $C^\ast$-algebra that is the
  pull-back of two unital, properly infinite $C^\ast$-algebras $A_1$ and
  $A_2$ along the $^\ast$-epimorphisms $\pi_1 \colon A_1 \to B$ and $\pi_2
  \colon A_2 \to B$:

\begin{displaymath}
\xymatrix{& A \ar[dl]_-{\varphi_1} \ar[dr]^-{\varphi_2} \\
A_1 \ar[dr]_-{\pi_1} && A_2 \ar[dl]^-{\pi_2}\\ 
& B &}
\end{displaymath}
Then $M_2(A)$ is properly infinite. Moreover, if $B$ is
$K_1$-injective, then $A$ itself is properly infinite. 
\end{proposition} 

\begin{proof}
Take unital embeddings $\sigma_i \colon \cT_3 \to A_i$ for $i=1,2$,
where $\cT_3$ is the Cuntz-Toeplitz algebra (defined earlier), 
and put
$$v = \sum_{j=1}^2 (\pi_1 \circ \sigma_1)(t_j) (\pi_2 \circ
\sigma_2)(t_j^*),$$
where $t_1,t_2,t_3$ are the canonical generators of $\cT_3$. Note that
$v$ is a partial isometry with $(\pi_1 
\circ \sigma_1)(t_j) = v (\pi_2 \circ \sigma_2)(t_j)$ for $j=1,2$.  
As $(\pi_1 \circ \sigma_1)(t_3t_3^*) \le 1-vv^*$ and 
$(\pi_2 \circ \sigma_2)(t_3t_3^*) \le 1-v^*v$,
Lemma~\ref{lm:unitary}~(i) yields a unitary $u \in B$ with $[u]=0$ in
$K_1(B)$ and with $(\pi_1
\circ \sigma_1)(t_j) = u (\pi_2 \circ \sigma_2)(t_j)$ for $j=1,2$.

If $B$ is $K_1$-injective, then $u$ belongs to $\cU^0(B)$,
whence $u$ lifts to a unitary $v \in A_2$. Define
$\widetilde{\sigma}_2 \colon \cT_2 \to 
A_2$ by $\widetilde{\sigma}_2(t_j) = v \sigma_2(t_j)$ for $j=1,2$
(observing that $t_1,t_2$ generate $\cT_2$). Then $\pi_1 \circ \sigma_1 =
\pi_2 \circ \widetilde{\sigma}_2$, which by the universal property of
the pull-back implies that $\sigma_1$ and $\widetilde{\sigma}_2$ lift
to a (necessarily unital) embedding $\sigma \colon \cT_2 \to A$, thus
forcing $A$ to be properly infinite.

In the general case (where $B$ is not necessarily
$K_1$-injective) $u$ may not lift to a unitary element in $A_2$, but
$\diag(u,u)$ does lift to a unitary element $v$ in $M_2(A_2)$  by
Lemma~\ref{lm:unitary}~(ii) (applied with $p = \diag(1,0)$). 
Define unital embeddings
$\widetilde{\sigma}_i \colon \cT_2 \to M_2(A_i)$, $i=1,2$,  by
$$\widetilde{\sigma}_1(t_j) = \begin{pmatrix} \sigma_1(t_j) & 0 \\0 &
  \sigma_1(t_j) \end{pmatrix}, \qquad  \widetilde{\sigma}_2(t_j) = v
\begin{pmatrix} \sigma_2(t_j) & 0 \\0 &
  \sigma_2(t_j) \end{pmatrix},$$
for $j=1,2$. As $(\pi_1 \otimes \Id_{M_2}) \circ \widetilde{\sigma}_1
= (\pi_2 \otimes \Id_{M_2}) \circ \widetilde{\sigma}_2$,  the 
unital embeddings $\widetilde{\sigma}_1$ and $\widetilde{\sigma}_2$ 
lift to a (necessarily unital) embedding of $\cT_2$
into $M_2(A)$, thus completing the proof.
\end{proof}

\begin{question} \label{q:pull-back} Is the pull-back of any two
  properly infinite unital \Cs s again properly infinite?
\end{question}

\noindent As mentioned in the introduction, one cannot in general 
conclude that $A$ is properly infinite if one knows that $M_n(A)$ 
is properly infinite for some $n \ge 2$. 

One obvious way of obtaining an answer to
Question~\ref{q:pull-back}, in the light of the last statement in
Proposition~\ref{prop:pull-back},  is to answer the question below in the
affirmative: 

\begin{question} \label{q:K_1-inj} Is every properly infinite unital
  \Cs{} $K_1$-injective?
\end{question}

\noindent
We shall see later, in Section~\ref{sec:equiv}, that the two questions
above in fact are equivalent.  

The lemma below, which shall be used several times in this paper,
shows that one can lift proper infiniteness from a fibre of a
$C(X)$-algebra to a whole neighborhood of that fibre.

\begin{lemma}\label{neighborhood}
Let $X$ be a compact Hausdorff space, let $A$ be a unital
$C(X)$-algebra, let $x \in X$, and suppose that the fibre $A_x$ is 
properly infinite. Then $A_{F}$ is properly infinite for some closed
neighborhood $F$ of $x$. 
\end{lemma}

\begin{proof}
Let $\{F_\lambda\}_{\lambda \in \Lambda}$ be a decreasing net of
closed neighborhoods of $x \in X$, fulfilling that
$\bigcap_{\lambda \in \Lambda} F_\lambda = \{x\}$, and set $I_\lambda =
C_0(X\setmin F_\lambda)A$. Then $\{I_\lambda\}_{\lambda \in \Lambda}$
is an increasing net of ideals in $A$, $A_{F_\lambda} = A/I_\lambda$, $I:=
\overline{\bigcup_{\lambda \in \Lambda} I_\lambda} = C_0(X \setmin
\{x\})$, and $A_x = A/I$. 

By the assumption that $A_x$ is properly infinite there is a unital
$^\ast$-homomorphism $\psi \colon \mathcal{T}_2 \to A_x$, and since
$\mathcal{T}_2$ is semi-projective there is a $\lambda_0 \in \Lambda$
and a unital $^\ast$-homomorphism $\varphi \colon \mathcal{T}_2 \to
A_{F_{\lambda_0}}$ making the diagram 
\begin{displaymath}
\xymatrix{& A_{F_{\lambda_0}} \ar[d]^-{\pi_x}\\
\mathcal{T}_2 \ar@{-->}[ur]^-{\varphi} \ar[r]_{\psi} & A_x}
\end{displaymath}
commutative. We can thus take $F$ to be $F_{\lambda_0}$. 
\end{proof}

\begin{theorem} \label{thm:C(X)-alg}
Let $A$ be a unital $C(X)$-algebra where $X$ is a compact
  Hausdorff space. If all fibres $A_x$, $x \in X$, are properly infinite,
  then some matrix algebra over $A$ is properly infinite.
\end{theorem}

\begin{proof} By Lemma \ref{neighborhood}, $X$ can be covered by
  finitely many closed sets $F_1, F_2, \dots, F_n$ such that $A_{F_j}$
  is properly infinite for each $j$. Put
  $G_j = F_1 \cup F_2 \cup \cdots \cup F_j$. For each $j=1, 2, \dots, n-1$ 
we have a pull-back diagram
\begin{displaymath}
\xymatrix@R-1pc@C-1pc{& A_{G_{j+1}} \ar[dl] \ar[dr] \\
A_{G_j} \ar[dr] && A_{F_{j+1}} \ar[dl]\\ 
&A_{G_j \cap F_{j+1}} &}
\end{displaymath}
We know that $M_{2^{j-1}}(A_{G_j})$ is properly infinite when $j=1$.
Proposition \ref{prop:pull-back} (applied to the diagram above
tensored with $M_{2^{j-1}}(\C)$) tells us that
$M_{2^{j}}(A_{G_{j+1}})$ is properly infinite if
$M_{2^{j-1}}(A_{G_{j}})$ is properly infinite. Hence $M_{2^{n-1}}(A)$ is
properly infinite.
\end{proof}

\begin{remark} 
Uffe Haagerup has suggested another way to prove Theorem \ref{thm:C(X)-alg}:  
If no matrix-algebra over $A$ is properly infinite, 
then there exists a bounded non-zero lower semi-continuous
$2$-quasi-trace on $A$, see \cite{DH} and \cite[page 327]{bh}, and
hence also an extremal $2$-quasi-trace. 
Now, if $A$ is also a $C(X)$-algebra for some compact Hausdorff space $X$, 
this implies that there is a bounded non-zero lower semi-continuous
$2$-quasitrace on $A_x$  
for (at least) one point $x\in X$ (see eg.\ \cite[Proposition~3.7]{hrw}). 
But then the fibre $A_x$ cannot be properly infinite. 
\end{remark} 

\begin{question} \label{q:C(X)} Is any unital $C(X)$-algebra $A$ properly
  infinite if all its fibres $A_x$, $x \in X$, are properly infinite?
\end{question} 

\noindent We shall show in Section~\ref{sec:equiv} that the question
above is equivalent to Question~\ref{q:pull-back} which again is
equivalent to Question~\ref{q:K_1-inj}. 

\section{Lower semi-continuous 
fields of properly infinite  
  $C^\ast$-algebras}  
Let us briefly discuss whether the results from Section~\ref{sec:C(X)} 
can be extended to lower
semi-continuous $C^\ast$-bundles $(A, \{\sigma_x\})$ over a
compact Hausdorff space~$X$.  
Recall that any such \textit{separable} lower semi-continuous
$C^\ast$-bundle admits
a faithful $C(X)$-linear representation on a Hilbert $C(X)$-module $E$
such that, for all $x\in X$, the fibre $\sigma_x(A)$ is isomorphic to
the induced image of $A$ in $\mathcal{L}(E_x)$, \cite{lsc}.  
Thus, the problem boils down to the following: 
Given a separable Hilbert $C(X)$-module $E$ with infinite dimensional
fibres $E_x$, such that
the unit $p$ of the $C^\ast$-algebra $\mathcal{L}_{C(X)}(E)$ of
bounded adjointable $C(X)$-linear operators acting on $E$ has a
properly infinite image in $\mathcal{L}(E_x)$ for all $x\in X$.  
Is the projection $p$ itself properly infinite in
$\mathcal{L}_{C(X)}(E)$?  

Dixmier and Douady proved that this is always the case if the
space $X$ has finite topological dimension, \cite{didou}.  
But it does not hold anymore in the infinite dimensional case, see
\cite[\S 16, Corollaire~1]{didou} and \cite{ro2}, not even if $X$ is 
contractible, \cite[Corollary~3.7]{BKbun}.  

\section{Two examples} \label{sec:examples}

\noindent We describe here two examples of continuous
fields; the first is over the interval and the second (which really is
a class of examples) is over the circle.  

\begin{example} \label{ex1}
Let $(\mathcal{O}_\infty \ast \mathcal{O}_\infty, (\iota_1, \iota_2))$ be
the universal unital free product of two copies of $\mathcal{O}_\infty$, 
and let $\mathcal{A}$ be
the unital sub-$C^\ast$-algebra of 
  $C([0,1],\mathcal{O}_\infty \ast \mathcal{O}_\infty)$ given by 
$$\mathcal{A} = \{f \in C([0,1], \mathcal{O}_\infty \ast
\mathcal{O}_\infty) : f(0) \in \iota_1(\mathcal{O}_\infty), \, f(1) \in
\iota_2(\mathcal{O}_\infty)\}.$$ 
Observe that $\mathcal{A}$ (in a canonical way) is a $C([0,1])$-algebra
with fibres
$$\mathcal{A}_t = \begin{cases} \iota_1(\cO_\infty), &
  t=0,\\  \mathcal{O}_\infty \ast \mathcal{O}_\infty, & 0 < t < 1,\\
\iota_2(\cO_\infty), & t=1
  \end{cases} \; \; \cong \; \; \begin{cases} \cO_\infty, &
  t=0,\\  \mathcal{O}_\infty \ast \mathcal{O}_\infty, & 0 < t < 1,\\
\cO_\infty, & t=1.
  \end{cases}$$
In particular, all fibres of $\mathcal{A}$ are properly infinite. 
\end{example}

\noindent One claim to fame of the example above is that the question
below is equivalent to Question~\ref{q:C(X)} above. Hence, to answer
Question~\ref{q:C(X)} in the affirmative (or in the negative) we need
only consider the case where $X= [0,1]$, and we need only worry about
this one particular $C([0,1])$-algebra (which of course is bad
enough!). 

\begin{question} \label{q:ex1}
Is the $C([0,1])$-algebra $\mathcal{A}$ from Example~\ref{ex1} above
properly infinite?
\end{question}

\noindent The three equivalent statements in the proposition below
will in Section~\ref{sec:equiv} be shown to be equivalent to
Question~\ref{q:ex1}. 

\begin{proposition} \label{prop:A} The following three statements
  concerning the  $C([0,1])$-algebra $\mathcal{A}$ and the 
\Cs{} $(\mathcal{O}_\infty \ast \mathcal{O}_\infty, (\iota_1, \iota_2))$
  defined above are equivalent:
\begin{enumerate}
\item $\cA$ contains a non-trivial projection (i.e., a projection
 other than $0$ and $1$).
\item There are non-zero projections $p,q \in \cO_\infty$ such that $p
  \ne 1$, $q \ne 1$, and $\iota_1(p) \sim_h \iota_2(q)$.
\item Let $s$ be any isometry in $\cO_\infty$. Then $\iota_1(ss^*)
  \sim_h \iota_2(ss^*)$ in $\cO_\infty \ast \cO_\infty$.
\end{enumerate}
\end{proposition}

\noindent We warn the reader that all three statements above could be
false.

\begin{proof} (i) $\Rightarrow$ (ii). Let $e$ be a non-trivial
  projection in $\cA$. Let $\pi_t \colon \cA \to \cA_t$, $t \in
  [0,1]$, denote the
  fibre map. As $\cA \subseteq C([0,1],\cO_\infty \ast \cO_\infty)$,
  the mapping $t \mapsto \pi_t(e) \in \cO_\infty \ast \cO_\infty$ is
  continuous, so in particular, $\pi_0(e) \sim_h \pi_1(e)$ in 
  $\cO_\infty \ast \cO_\infty$. The mappings $\iota_1$ and
  $\iota_2$ are injective, so there are projections $p,q \in \cO_\infty$
  such that $\pi_0(e) = \iota_1(p)$ and $\pi_1(e) = \iota_2(q)$. The
  projections $p$ and $q$ are non-zero because the mapping $t \mapsto
  \|\pi_t(e)\|$ is continuous and not constant equal to 0. Similarly,
  and $1-p$ and $1-q$ are non-zero because $1-e$ is non-zero. 

(ii) $\Rightarrow$ (iii). Take non-trivial projections $p,q \in
\cO_\infty$ such that $\iota_1(p) \sim_h \iota_2(q)$. Take a unitary
$v$ in $\cU^0(\cO_\infty \ast \cO_\infty)$ with $\iota_2(q) =
v\iota_1(p)v^*$. Let $s \in \cO_\infty$ be an isometry. If $s$ is
unitary, then $\iota_1(ss^*) = 1 = \iota_2(ss^*)$ and there is nothing
to prove. Suppose that $s$ is non-unitary. Then $ss^*$ is homotopic to
a subprojection $p_0$ of $p$ and to a subprojection $q_0$ of $q$ 
(use that $p$ and $q$ are properly infinite and full, then
Lemma~\ref{lm:unitary}~(i), and last the fact 
that the unitary group of $\cO_\infty$ is connected). Hence
$\iota_1(ss^*) \sim_h \iota_1(p_0) \sim_h v\iota_1(p_0)v^*$ and
$\iota_2(ss^*) \sim_h \iota_2(q_0)$, so we
need only show that $v \iota_1(p_0)v^* \sim_h \iota_2(q_0)$. But this
follows from Proposition~\ref{prop:propinf2} with $r= 1
-\iota_2(q) = \iota_2(1-q)$, 
as we note that $p_0 \sim 1 \sim q_0$ in $\cO_\infty$,
whence $$\iota_2(q_0) \sim \iota_2(1) = 1 = \iota_1(1) \sim
\iota_1(p_0) \sim v\iota_1(p_0)v^*.$$

(iii) $\Rightarrow$ (i). Take a non-unitary isometry $s \in
\cO_\infty$. Then $\iota_1(ss^*) \sim_h \iota_2(ss^*)$, and so there
is a continuous function $e \colon [0,1] \to \cO_\infty \ast
\cO_\infty$ such that $e(t)$ is a projection for all $t \in [0,1]$,
$e(0) = \iota_1(ss^*)$ and $e(1) = \iota_2(ss^*)$. But then $e$ is a
non-trivial projection in $\cA$. 
\end{proof}

\noindent It follows from Theorem~\ref{thm:C(X)-alg} that some matrix
algebra over $\cA$ (from Example~\ref{ex1}) 
is properly infinite. We can sharpen that statement as follows:

\begin{proposition} \label{prop:M_2(A)}
$M_2(\cA)$ is properly infinite; and if $\cO_\infty \ast \cO_\infty$
is $K_1$-injective, then $\cA$ itself is properly infinite.
\end{proposition}

\noindent It follows from Theorem~\ref{equiv conditions} below that
$\cA$ is properly infinite if \emph{and only if} 
$\cO_\infty \ast \cO_\infty$ is $K_1$-injective. 

\begin{proof} We have a pull-back diagram
\begin{equation*} 
\xymatrix{& \cA \ar[ld] \ar[rd] & \\ \cA_{[0,\frac{1}{2}]} \ar[rd]_{\pi_{1/2}}
  && \cA_{[\frac{1}{2},1]} \ar[ld]^{\pi_{1/2}} \\ & \cO_\infty \ast \cO_\infty
  &}
\end{equation*}
One can unitally embed $\cO_\infty$ into $\cA_{[0,\frac{1}{2}]}$ via
$\iota_1$, so $\cA_{[0,\frac{1}{2}]}$ is properly infinite, and a
similar argument shows that $\cA_{[\frac{1}{2},1]}$ is properly
infinite. The two statements now follow from Proposition~\ref{prop:pull-back}.
\end{proof}

\noindent 
The example below, which will be the focus of the rest of this section, and
in parts also of Section \ref{sec:equiv}, is inspired by arguments
from Rieffel's paper \cite{Rieffel}. 

\begin{example} \label{ex2} 
Let $A$ be a unital \Cs, and let $v$ be a unitary element in $A$ such that 
$$\begin{pmatrix} v & 0 \\ 0 & 1 \end{pmatrix} \sim_h \begin{pmatrix} 1 & 0
  \\ 0 & 1 \end{pmatrix}\quad \text{in} \quad \mathcal{U}_2(A).$$
Let $t \mapsto u_t$ be a continuous path of unitaries in $\mathcal{U}_2(A)$
such that $u_0 = 1$ and $u_1 = \diag(v,1)$. Put
$$p(t) = u_t  \begin{pmatrix} 1 & 0 \\ 0 & 0 \end{pmatrix} u_t^* \in
M_2(A),$$
and note that $p(0)=p(1)$. Identifying, for each \Cs{} $D$, $C(\T,D)$ with
the algebra of all continuous functions $f \colon [0,1] \to D$ such that
$f(1)=f(0)$, we see that $p$ belongs to $C(\T,M_2(A))$. Put
$$\mathcal{B} = p C(\T,M_2(A)) p,$$
and note that $\mathcal{B}$ is a unital (sub-trivial) $C(\T)$-algebra,
being a corner of the trivial $C(\T)$-algebra
$C(\T,M_2(A))$. The fibres of $\mathcal{B}$ are
$$\mathcal{B}_t = p(t)M_2(A)p(t) \cong A$$
for all $t \in \T$. 

Summing up, for each unital \Cs{} $A$, for each unitary $v$ in $A$ for
which $\diag(v,1) \sim_h 1$ in $\cU_2(A)$, and for each path $t
\mapsto u_t \in \cU_2(A)$ implementing this homotopy we get a
$C(\T)$-algebra $\cB$ with fibres $\cB_t \cong A$. We shall
investigate this class of $C(\T)$-algebras below.
\end{example}

\begin{lemma}\label{stable eq} In the notation of Example~\ref{ex2},
  $$\begin{pmatrix}1 & 0 \\ 0 & 1 \end{pmatrix} -p \; \sim \;  
\begin{pmatrix} 0 & 0 \\ 0 & 1 \end{pmatrix} \quad \text{in} \quad 
C(\T,M_2(A)).$$ 
In particular, $p$ is stably equivalent to $\diag(1,0)$.  
\end{lemma}

\begin{proof}
Put $$v_t = u_t\left(\begin{array}{cc}
0 & 0\\
0 & 1
\end{array}\right), \quad t \in [0,1].$$ Then
$$v_0 = u_0\left(\begin{array}{cc}
0 & 0\\
0 & 1
\end{array}\right) = \left(\begin{array}{cc}
0 & 0\\
0 & 1
\end{array}\right), \quad v_1 = u_1\left(\begin{array}{cc}
0 & 0\\
0 & 1
\end{array}\right) = \left(\begin{array}{cc}
v & 0\\
0 & 1
\end{array}\right) \left(\begin{array}{cc}
0 & 0\\
0 & 1
\end{array}\right)= \left(\begin{array}{cc}
0 & 0\\
0 & 1
\end{array}\right),$$
so $v$ belongs to $C(\T, M_2(A))$. It is easy to see that $v_t^*v_t =
\diag(0,1)$ and $v_t v_t^* =
1-p(t)$, and so the lemma is proved. 
\end{proof}

\begin{proposition}\label{prop1}  Let $A$, $v \in \cU(A)$, and 
$\mathcal{B}$ be as in Example~\ref{ex2}.
Conditions {\rm{(i)}} and {\rm{(ii)}} below are equivalent
  for any unital \Cs{} $A$, and all three conditions are equivalent
if $A$ in addition is assumed to be properly infinite.
\begin{enumerate}
\item $v \sim_h 1$ in $\mathcal{U}(A)$.
\item $p \sim \diag(1_A,0)$ in $C(\T,M_2(A))$.
\item The $C(\T)$-algebra $\mathcal{B}$ is properly infinite. 
\end{enumerate}
\end{proposition}

\begin{proof}
(ii) $\Rightarrow$ (i).
Suppose that $p \sim  \diag(1,0)$ in $C(\T,M_2(A))$. Then there is
a $w \in C(\T, M_2(A))$ such that  
$$w_tw_t^* = \left ( \begin{array}{cc}
1 & 0\\
0 & 0
\end{array} \right)\quad \text{and} \quad
w_t^\ast w_t  = p_t$$
for all $t \in [0,1]$ and $w_1=w_0$ (as we identify $C(\T,M_2(A))$
with the set of continuous functions $f \colon [0,1] \to M_2(A)$ with
$f(1)=f(0)$). Upon replacing $w_t$ with $w_0^*w_t$ we can assume
that $w_1=w_0=\diag(1,0)$.  Now, with $t \mapsto u_t$ as in
Example~\ref{ex2}, 
$$w_t u_t \left( \begin{array}{cc}
1 & 0\\
0 & 0
\end{array}\right) = \begin{pmatrix} a_t & 0 \\ 0 & 0 \end{pmatrix},$$
where $t \mapsto a_t$ is a continuous path of unitaries in $A$. Because
$u_0= \diag(1,1)$ and $u_1= \diag(v,1)$ we see that $a_0 = 1$ and $a_1
= v$, whence $v \sim_h 1$ in $\cU(A)$. 

(i) $\Rightarrow$ (ii). Suppose conversely that $v \sim_h 1$ in
$\mathcal{U}(A)$. Then we can find a continuous path $t \mapsto v_t
\in \cU(A)$, $t \in [1-\ep,1]$, such that $v_{1-\varepsilon} = v$ and $v_1
= 1$ for an $\varepsilon >0$ (to be determined below). Again with
$t\mapsto u_t$ as in Example~\ref{ex2}, define
  \begin{displaymath}
\widetilde{u}_t = \begin{cases}
u_{(1-\varepsilon)^{-1} t}, & 0 \le t \le 1-\ep,\\
\diag(v_t,1), & 1-\ep \le t \le 1.
\end{cases}
\end{displaymath}
Then $t \mapsto \widetilde{u}_t$ is a continuous path of unitaries in
$\mathcal{U}_2(A)$ such that $\widetilde{u}_{1-\varepsilon} = u_1 =
\diag(v,1)$ and $\widetilde{u}_0=\widetilde{u}_1 = 1$. It follows that
$\widetilde{u}$ 
belongs to $C(\T,M_2(A))$. Provided that $\varepsilon>0$
is chosen small enough we obtain the following inequality: 
$$\left\|\widetilde{u}_t\left(\begin{array}{cc}
1 & 0\\
0 & 0
\end{array}\right)\widetilde{u}_t ^\ast - p(t)\right\| =
\left\|\widetilde{u}_t \left(\begin{array}{cc} 
1 & 0\\
0 & 0
\end{array}\right)\widetilde{u}_t ^\ast - u_t\left(\begin{array}{cc}
1 & 0\\
0 & 0
\end{array}\right)u_t^\ast\right\| <1$$
for all $t \in [0,1]$, whence $p 
\sim \widetilde{u} \, \diag(1,0) \,
\widetilde{u}^* \sim \diag(1,0)$ as desired. 

(iii) $\Rightarrow$ (ii). Suppose that $\mathcal{B}$ is properly
infinite. From Lemma \ref{stable eq} we know that $[p] =
\left[\diag(1_A,0)\right]$ in $K_0(C(\T,A))$. Because $\cB$ and $A$
are properly infinite, it follows that $p$ and $\diag(1_A,0)$ are
properly infinite (and full) projections, and hence they are 
equivalent by Proposition~\ref{prop:propinf1}~(i).

(ii) $\Rightarrow$ (iii). Since $A$ is properly infinite,
$\diag(1_A,0)$ and hence $p$ (being equivalent to $\diag(1_A,0)$) 
are properly infinite (and full) 
projections, whence $\cB$ is properly infinite.  
\end{proof}

\noindent We will now use (the ideas behind)
Lemma~\ref{stable eq} and Proposition
\ref{prop1} to prove the following general statement about
$C^\ast$-algebras. 

\begin{corollary}\label{kappa2}
Let $A$ be a unital $C^\ast$-algebra such that $C(\T,A)$ has the
cancellation property. Then $A$ is $K_1$-injective. 
\end{corollary}

\begin{proof} It suffices to show that the natural maps
  $\mathcal{U}_{n-1}(A)/\mathcal{U}_{n-1}^0(A) \to 
\mathcal{U}_n(A)/\mathcal{U}_n^0(A)$ are injective for all $n \geq 2$.
Let $v \in \mathcal{U}_{n-1}(A)$ be
such that $\diag(v,1_A) \in \mathcal{U}_n^0(A)$ and
find a continuous path of unitaries $t \mapsto u_t$ in
$\mathcal{U}_n(A)$ such that  
$$u_0 = 1_{M_n(A)} = \left(\begin{array}{cc}
1_{M_{n-1}(A)} & 0\\
0 & 1_A
\end{array}\right)\quad \text{and} \quad u_1 = \left(\begin{array}{cc}
v & 0\\
0 & 1_A
\end{array}\right).$$ Put
$$p_t = u_t \left(\begin{array}{cc}
1_{M_{n-1}(A)} & 0\\
0 & 0 \end{array}\right) u_t^\ast, \quad t \in [0,1],$$ 
and note that $p_0=p_1$ so that 
$p$ defines a projection in $C(\T, M_n(A))$. Repeating
the proof of Lemma~\ref{stable eq} we find that $1_{M_n(A)} - p \sim
\diag(0,1_A)$ in $C(\T,M_n(A))$, whence $p \sim
\diag(1_{M_{n-1}(A)},0)$ by the cancellation property of $C(\T,A)$,
where we identify projections in $M_n(A)$ with
constant projections in $C(\T,M_n(A))$. The
arguments going into the proof of Proposition~\ref{prop1} show that $v
\sim_h 1_{M_{n-1}(A)}$ in $\mathcal{U}_{n-1}(A)$ if (and only if) $p
\sim \diag (1_{M_{n-1}(A)},0)$. Hence $v$ belongs to $\cU_{n-1}^0(A)$
as desired.
\end{proof}

\section{$K_1$-injectivity of properly infinite \Cs s} \label{sec:equiv}

In this section we prove our main result that relate $K_1$-injectivity
of arbitrary unital properly infinite \Cs s to proper infiniteness of
$C(X)$-algebras and pull-back \Cs s. More specifically we shall show
that  Question~\ref{q:K_1-inj}, Question~\ref{q:C(X)},
Question~\ref{q:pull-back}, and Question~\ref{q:ex1} are equivalent.  

First we reformulate in two different ways
the question if a given properly infinite unital
\Cs{} is $K_1$-injective.

\begin{proposition}\label{prop:equiv}
The following conditions are equivalent for any unital properly infinite
\Cs{} $A$:
\begin{enumerate}
\item $A$ is $K_1$-injective.
\item Let $p$, $q$ be projections in $A$ such that $p \sim q$ and $p$,
  $q$, $1-p$, $1-q$ are properly infinite and full. Then $p \sim_h q$. 
\item Let $p$ and $q$ be properly infinite, full projections in
  $A$. There exist properly infinite, full projections $p_0, q_0 \in A$
  such that $p_0 \leq p$, $q_0 \leq q$, and $p_0 \sim_h q_0$. 
\end{enumerate}
\end{proposition}

\begin{proof}
(i) $\Rightarrow$ (ii). Let $p,q$ be properly infinite, full
projections in $A$ with $p \sim q$ such that $1-p, 1-q$ are properly
infinite and full. Then by  Lemma~\ref{lm:unitary}~(i) 
there is a unitary $v \in A$ such that $vpv^\ast =q$ and $[v]=0$ in
$K_1(A)$. By the assumption in (i), $v \in \cU^0(A)$, whence $p \sim_h
q$. 

(ii) $\Rightarrow$ (i). Let $u\in \mathcal{U}(A)$ be such that $[u] =
0$ in $K_1(A)$. Take, as we can, a projection $p$ in $A$ such that $p$
and $1-p$ are properly infinite and full. Set $q= upu^\ast$. Then $p \sim_h
q$ by (ii), and so there exists a unitary $v \in 
\mathcal{U}^0(A)$ with  $p= vqv^\ast$. It follows that 
$$pvu =vqv^*vu = v(upu^*) v^* vu = vup.$$
Therefore $vu \in \mathcal{U}^0(A)$ by
Lemma~\ref{lm:unitary}~(ii), 
which in turn implies that $u \in \cU^0(A)$. 

(ii) $\Rightarrow$ (iii). Let $p,q$ be properly infinite and full
projections in $A$. There exist mutually orthogonal projections $e_1,
f_1$ such that $e_1 \leq p$, $f_1 \leq p$ and $e_1 \sim p \sim f_1$,
and mutually orthogonal projections $e_2, f_2$ such that $e_2 \leq q$,
$f_2 \leq q$ and $e_2 \sim q \sim f_2$. Being equivalent to either $p$ or
$q$, the projections $e_1,e_2,f_1$ and $f_2$ are properly infinite and
full. There are properly infinite, full
projections $p_0 \le e_1$ and $q_0 \le e_2$ such that $[p_0]=[q_0] =
0$ in $K_0(A)$  and $p_0 \sim q_0$ (cf.\ Proposition~\ref{prop:propinf1}). As
$f_1 \le 1-p_0$ and $f_2 \le 1-q_0$, we see that $1-p_0$ and
$1-q_0$ are properly infinite and full, and so we get $p_0 \sim_h q_0$
by (ii). 

(iii) $\Rightarrow$ (ii). Let $p,q$ be equivalent properly infinite, full
projections in $A$ such that 
$1-p, 1-q$ are properly infinite and full.
From (iii) we get properly infinite and full
projections $p_0 \leq p$, $q_0 \leq q$ which satisfy $p_0 \sim_h
q_0$. Thus there is a unitary $v \in \mathcal{U}_0(A)$ such that
$vp_0v^\ast = q_0$. Upon
replacing $p$ by $vpv^*$ (as we may do because $p \sim_h vpv^*$) we
can assume that $q_0 \le p$ and $q_0 \le q$. Now, $q_0$ is orthogonal
to $1-p$ and to $1-q$, and so $1-p \sim_h 1-q$ by
Proposition~\ref{prop:propinf2}, whence $p \sim_h q$. 
\end{proof}

\begin{proposition} \label{prop:K1inj}
Let $A$ be a unital properly infinite \Cs. The following conditions
are equivalent:
\begin{enumerate}
\item $A$ is $K_1$-injective, ie., the natural map
  $\cU(A)/\cU^0(A) \to K_1(A)$ is injective.
\item The natural map $\cU(A)/\cU^0(A) \to \cU_2(A)/\cU_2^0(A)$ is
  injective. 
\item The natural maps $\cU_n(A)/\cU_n^0(A) \to K_1(A)$ are
  injective for each natural number $n$. 
\end{enumerate}
\end{proposition}

\begin{proof} (i) $\Rightarrow$ (ii) holds because the map
  $\cU(A)/\cU^0(A) \to K_1(A)$ factors through the map
  $\cU(A)/\cU^0(A) \to \cU_2(A)/\cU_2^0(A)$.

(ii) $\Rightarrow$ (i). Take $u \in \cU(A)$ and suppose that $[u]=0$ in
$K_1(A)$. Then $\diag(u,1_A) \in \cU_2^0(A)$ by
Lemma~\ref{lm:unitary}~(ii) (with $p = \diag(1_A,0)$). Hence $u
\in \cU_0(A)$ by injectivity of the map $\cU(A)/\cU^0(A) \to
\cU_2(A)/\cU_2^0(A)$. 

(i) $\Rightarrow$ (iii). Let $n \ge 1$ be given and consider the natural maps
$$\cU(A)/\cU^0(A) \to \cU_n(A)/\cU_n^0(A) \to K_1(A).$$
The first map is onto, as proved by Cuntz in \cite{Cu}, see also
\cite[Exercise 8.9]{LLR}, 
and the composition of the two maps is injective by assumption, hence
the second map is injective.

(iii) $\Rightarrow$ (i) is trivial.
\end{proof}

\noindent We give below another application of $K_1$-injectivity for
properly infinite \Cs s. First we need a lemma:

\begin{lemma} \label{lm:homotopy}
Let $A$ be a unital, properly infinite \Cs, 
and let $\varphi, \psi \colon \cO_\infty \to
A$ be unital embeddings. Then $\psi$ is homotopic to a unital
embedding $\psi' \colon \cO_\infty \to A$ for which there is a
unitary $u \in A$ with $[u]=0$ in $K_1(A)$ and for which $\psi'(s_j) = u
\varphi(s_j)$ for all $j$ (where $s_1,s_2, \dots$ are the canonical
generators of $\cO_\infty$).
\end{lemma}

\begin{proof} For each $n$ set
$$v_n = \sum_{j=1}^n \psi(s_j) \varphi(s_j)^* \in A, \qquad e_n =
\sum_{j=1}^n s_js_j^* \in \cO_\infty.$$
Then $v_n$ is a partial isometry in $A$ with $v_nv_n^* = \psi(e_n)$,
$v_n^*v_n = \varphi(e_n)$, and $\psi(s_j) = v_n\varphi(s_j)$ for
$j=1, 2, \dots, n$. Since $1-e_n$ is full and properly infinite it
follows from Lemma~\ref{lm:unitary} that each $v_n$ extends to a unitary
$u_n \in A$ with $[u_n]=0$ in $K_1(A)$. In particular, $\psi(s_j) =
u_n\varphi(s_j)$ for $j=1,2, \dots, n$. 

We proceed to show that $n \mapsto u_n$ extends to a continuous path
of unitaries $t \mapsto u_t$, for $t \in [2,\infty)$, such that $u_t
\varphi(e_n) = u_n \varphi(e_n)$ for $t \ge n+1$. Fix $n \ge 2$. To
this end it
suffices to show that we can find a continuous path $t \mapsto z_t$,
$t \in [0,1]$, of unitaries in $A$ 
such that $z_0 = 1$, $z_1 = u_{n}^*u_{n+1}$, and
$z_t\varphi(e_{n-1}) = \varphi(e_{n-1})$ (as we then can set $u_t$ to
be $u_nz_{t-n}$ for $t \in [n,n+1]$). 

Observe that 
$$u_{n+1}\varphi(e_n) = v_{n+1}\varphi(e_n) = v_n = u_n\varphi(e_n).$$
Set $A_{0} = (1-\varphi(e_{n-1}))A(1-\varphi(e_{n-1}))$, and set
$y = u_n^*u_{n+1}(1- \varphi(e_{n-1}))$. Then $y$ is a unitary
element in $A_{0}$ and $[y] = 0$ in $K_1(A_{0})$. Moreover,
$y$ commutes  with the properly infinite full projection
$\varphi(e_n)-\varphi(e_{n-1}) \in A_{0}$. We can therefore use
Lemma~\ref{lm:unitary} to find a continuous path $t \mapsto y_t$ of
unitaries in $A_{0}$ such that
$y_0 = 1_{A_{0}} = 1- \varphi(e_{n-1})$ and $y_1 =y$. The
continuous path $t \mapsto z_t = y_t + \varphi(e_{n-1})$ is then as desired.

For each $t \ge 2$ let $\psi_t \colon \cO_\infty \to A$ be the \sh{}
given by $\psi_t(s_j) = u_t \varphi(s_j)$. Then $\psi_t(s_j) =
\psi(s_j)$ for all $t \ge j+1$, and so it follows that
$$\lim_{t \to \infty} \psi_t(x) = \psi(x)$$
for all $x \in \cO_\infty$. Hence $\psi_2$ is homotopic to $\psi$, and
so we can take $\psi'$ to be $\psi_2$.
\end{proof}

\begin{proposition} \label{prop:homotopy1}
Any two unital \sh s from $\cO_\infty$ into a unital $K_1$-injective
(properly infinite) \Cs{} are homotopic. 
\end{proposition}

\begin{proof}
In the light of  Lemma~\ref{lm:homotopy} it suffices to show that if
$\varphi, \psi \colon \cO_\infty \to A$ are unital \sh s such that,
for some unitary $u \in A$ with $[u]=0$ in $K_1(A)$, $\psi(s_j) =
u\varphi(s_j)$ for all $j$, then $\psi \sim_h \varphi$. By assumption,
$u \sim_h 1$, so there is a continuous path $t \mapsto u_t$ of
unitaries in $A$ such that $u_0=1$ and $u_1=u$. Letting $\varphi_t
\colon \cO_\infty \to A$ be the \sh{} given by $\varphi_t(s_j) =
u_t\varphi(s_j)$ for all $j$, we get $t \mapsto \varphi_t$ is a
continuous path of \sh s connecting $\varphi_0=\varphi$ to
$\varphi_1=\psi$. 
\end{proof}

\noindent Our main theorem below, which in particular implies that  
Question \ref{q:K_1-inj}, Question~\ref{q:C(X)},
Question~\ref{q:pull-back} and Question~\ref{q:ex1} all are
equivalent, also give a special converse to
Proposition~\ref{prop:homotopy1}: Indeed, with $\iota_1, \iota_2
\colon \cO_\infty \to \cO_\infty \ast \cO_\infty$ the two canonical
inclusions, if $\iota_1 \sim_h \iota_2$, then condition (iv) below
holds, whence $\cO_\infty \ast \cO_\infty$ is $K_1$-injective, which
again implies that all unital properly infinite \Cs s are
$K_1$-injective. Below we retain the convention that $\cO_\infty \ast
\cO_\infty$ is the universal \emph{unital} free product of two copies
of $\cO_\infty$ and that $\iota_1$ and $\iota_2$ are the two natural
inclusions of $\cO_\infty$ into $\cO_\infty \ast \cO_\infty$. 

\begin{theorem} \label{equiv conditions}
The following statements are equivalent:
\begin{enumerate}
\item Every unital, properly infinite $C^\ast$-algebra is $K_1$-injective.
\item For every compact Hausdorff space $X$, every unital
  $C(X)$-algebra $A$, for which $A_x$ is properly infinite for all $x \in
  X$, is properly infinite. 
\item Every unital $C^\ast$-algebra $A$, that is the pull-back of
  two unital, properly infinite $C^\ast$-algebras $A_1$ and $A_2$
  along $^\ast$-epimorphisms $\pi_1 \colon A_1 \to B$,
  $\pi_2\colon A_2 \to B$:

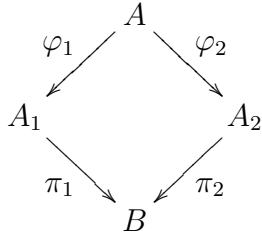
\begin{displaymath}
\xymatrix{& A \ar[dl]_-{\varphi_1} \ar[dr]^-{\varphi_2} \\
A_1 \ar[dr]_-{\pi_1} && A_2 \ar[dl]^-{\pi_2}\\ 
& B &}
\end{displaymath}
is properly infinite.
\item There exist non-zero projections $p,q \in
  \mathcal{O}_\infty$ such that $p \ne 1$, $q \ne 1$, and $\iota_1(p) \sim_h
  \iota_2(p)$ in $\mathcal{O}_\infty \ast \mathcal{O}_\infty$.
\item The specific $C([0,1])$-algebra $\mathcal{A}$ considered in
  Example \ref{ex1} (and whose fibres are properly infinite) 
is properly infinite. 
\item $\cO_\infty \ast \cO_\infty$ is $K_1$-injective. 
\end{enumerate}
\end{theorem}

\noindent Note that statement (i) is reformulated in
Propositions~\ref{prop:equiv}, \ref{prop:K1inj}, and
\ref{prop:homotopy1}; and that 
statement (iv) is reformulated in
Proposition~\ref{prop:A}. We warn the reader that all these statements
may turn out to be false (in which case, of course, there will be
counterexamples to all of them). 

\begin{proof} (i) $\Rightarrow$ (iii) follows from
  Proposition~\ref{prop:pull-back}. 

(iii) $\Rightarrow$ (ii). This follows from Lemma~\ref{neighborhood}
as in the proof of 
Theorem~\ref{thm:C(X)-alg}, except that one does not need to pass
to matrix algebras.   

(ii) $\Rightarrow$ (i). Suppose that $A$ is unital and properly
infinite. Take a unitary $v \in \cU(A)$ such that $\diag(v,1) \in
\cU^0_2(A)$. Let $\mathcal{B}$ be the $C(\T)$-algebra constructed in
Example~\ref{ex2} from $A$, $v$, and a path of unitaries $t \mapsto u_t$
connecting $1_{M_2(A)}$ to $\diag(v,1)$. Then $\mathcal{B}_t \cong A$
for all $t \in \T$, so all fibres of $\mathcal{B}$ are properly
infinite. Assuming (ii), we can conclude that $\mathcal{B}$ is
properly infinite.  
Proposition~\ref{prop1} then yields that $v \in \cU^0(A)$. It follows
that the natural map $\cU(A)/\cU_0(A) \to \cU_2(A)/\cU_2^0(A)$
is injective, whence $A$ is $K_1$-injective by
Proposition~\ref{prop:K1inj}. 

(ii) $\Rightarrow$ (v) is trivial (because $\mathcal{A}$ is a
$C([0,1])$-algebra with properly infinite fibres). 

(v) $\Rightarrow$ (iv) follows from Proposition~\ref{prop:A}. 

(iv) $\Rightarrow$ (i). We show that Condition~(iii) of
Proposition~\ref{prop:A} implies Condition~(iii) of
Proposition~\ref{prop:equiv}.
 
Let $A$ be a properly infinite $C^\ast$-algebra and let $p,q$ be
properly infinite, full projections in $A$. Then there exist
(properly infinite, full) 
projections $p_0 \leq p$ and $q_0 \leq q$ such that $p_0 \sim 1 \sim
q_0$ and such that $1-p_0$ and $1-q_0$ are properly infinite and full,
cf.\ Propositions~\ref{prop:propinf1}. Take
isometries $t_1, r_1 \in A$ with $t_1 t_1^* = p_0$ and $r_1r_1^* =
q_0$; use the fact that $1 \precsim 1-p_0$ and $1 \precsim 1-q_0$
to find sequences of isometries $t_2,t_3,t_4, \dots$ and $r_2,r_3,r_4,
\dots$ in $A$ such that each of the two sequences
$\{t_jt_j^*\}_{j=1}^\infty$ and $\{r_jr_j^*\}_{j=1}^\infty$ consist
of pairwise orthogonal projections. 

By the universal property of $\cO_\infty$ there are unital
$^\ast$-homomorphisms 
$\varphi_j \colon \mathcal{O}_\infty \to A$, $j=1,2$, such that
$\varphi_1(s_j) = t_j$ and $\varphi_2(s_j)=r_j$, where $s_1,s_2,s_3,
\dots$ are the canonical generators of $\cO_\infty$. In particular, 
$$\varphi_1(s_1s_1^\ast) = p_0 \quad \text{and} \quad
\varphi_2(s_1s_1^\ast) = q_0.$$ 
By the property of the universal unital free products of
$C^\ast$-algebras, there is a unique unital $^\ast$-homomorphism
$\varphi \colon \mathcal{O}_\infty \ast \mathcal{O}_\infty \to A$ 
making the diagram 
\begin{displaymath}
\xymatrix{& \mathcal{O}_\infty \ast \mathcal{O}_\infty \ar[dd]_-{\varphi}\\
\mathcal{O}_\infty \ar[dr]_-{\varphi_1} \ar[ur]^-{\iota_1} &&
\mathcal{O}_\infty \ar[dl]^-{\varphi_2} \ar[ul]_-{\iota_2}\\  
& A &}
\end{displaymath}
commutative. It follows that $p_0 = \varphi(\iota_1(s_1s_1^\ast))$ and
$q_0 = \varphi(\iota_2(s_1s_1^\ast))$. By Condition~(iii) of
Proposition~\ref{prop:A}, $\iota_1(s_1s_1^\ast)
\sim_h \iota_2(s_1s_1^\ast)$ in $\cO_\infty \ast \cO_\infty$, whence
 $p_0 \sim_h q_0$ as desired.

(i) $\Rightarrow$ (vi) is trivial.

(vi) $\Rightarrow$ (v) follows from Proposition~\ref{prop:M_2(A)}.
\end{proof}

\section{Concluding remarks} 
\label{sec:moreexamples}

\noindent We do not know if all unital properly infinite \Cs s are
$K_1$-injective, but we observe that 
$K_1$-injectivity is assured in the presence of certain central sequences:

\begin{proposition} \label{prop:centralseq}
Let $A$ be a unital properly infinite \Cs s that contains 
an asymptotically central sequence $\{p_n\}_{n=1}^\infty$, where $p_n$
and $1-p_n$  are properly infinite, full projections for all $n$. Then $A$
is $K_1$-injective
\end{proposition}

\begin{proof} This follows immediately from Lemma~\ref{lm:unitary}~(ii).
\end{proof}

\noindent It remains open if arbitrary $C(X)$-algebras with properly
infinite fibres must be properly infinite. If this fails, then we
already have a counterexample of the form $\cB = pC(\T,A)p$, cf.\
Example~\ref{ex2}, for some unital properly infinite \Cs{} $A$ and
for some projection $p \in C(\T,A)$. (The \Cs{} $\cB$ is a
$C(\T)$-algebra with fibres $\cB_t \cong A$.) 

On the other hand, any trivial $C(X)$-algebra $C(X,D)$ with constant
fibre $D$ is clearly properly infinite if its fibre(s) $D$ is unital
and properly infinite (because $C(X,D) \cong C(X) \otimes D$). We
extend this observation in the following easy:

\begin{proposition} \label{prop:trivial}
Let $X$ be a compact Hausdorff space, let $p \in C(X,D)$ be a
projection, and consider the sub-trivial $C(X)$-algebra $pC(X,D)p$
whose fibre at $x$ is equal to $p(x)Dp(x)$. 

If $p$ is Murray-von Neumann equivalent to a constant projection $q$, then
$pC(X,D)p$ is $C(X)$-isomorphic to the trivial $C(X)$-algebra
$C(X,D_0)$, where $D_0 = qDq$. In this case, $pC(X,D)p$ is properly
infinite if and only if $D_0$ is properly infinite.

In particular, if $X$ is contractible, then $pC(X,D)p$ is $C(X)$-isomorphic
to a trivial $C(X)$-algebra for any projection $p \in C(X,D)$ and for
any \Cs{} $D$.
\end{proposition}

\begin{proof} Suppose that $p=v^*v$ and $q = vv^*$ for some partial
  isometry $v \in C(X,D)$. The map $f \mapsto vfv^*$ defines a
  $C(X)$-isomorphism from $pC(X,D)p$ onto $qC(X,D)q$, and $qC(X,D)q =
  C(X,D_0)$. 

If $X$ is contractible, then any projection $p \in C(X,D)$ is
homotopic, and hence 
equivalent, to the constant projection $x \mapsto p(x_0)$ for any
fixed $x_0 \in X$. 
\end{proof}

\begin{remark} \label{rem:subtrivial}
One can elaborate a little more on the construction considered
above. Take a unital \Cs{} $D$ such that for some natural number $n
\ge 2$, $M_n(D)$ is properly infinite, but $M_{n-1}(D)$ is not
properly infinite (see \cite{Ror:sums} or \cite{ro2} for such examples). Take
any space $X$, preferably one with highly non-trivial topology, eg.\
$X=S^n$, and take, for some $k \ge n$, a sufficiently non-trivial 
$n$-dimensional projection $p$ in $C(X,M_k(D))$ such that $p(x)$ is
equivalent to the trivial $n$ dimensional projection $1_{M_n(D)}$ for
  all $x$ (if $X$ is connected we need only assume that this holds for
  one $x \in X$). The $C(X)$-algebra 
$$\cA = p \, C(X,M_k(D)) \, p,$$
then has properly infinite fibres $\cA_x = p(x)Dp(x) \cong M_n(D)$. 
Is $A$ always properly infinite?
We guess that a possible counterexample to the questions posed in this paper 
could be of this form (for suitable $D$, $X$, and $p$). 
\end{remark}

\noindent Let us end this paper by remarking that the answer to 
Question~\ref{q:C(X)}, which asks if any $C(X)$-algebra with properly
infinite fibres is itself properly infinite, does not depend (very
much) on $X$. If it fails, then it fails already for $X = [0,1]$ (cf.\
Theorem~\ref{equiv conditions}), and $[0,1]$ is a contractible space
of low dimension. However, if we make the dimension of $X$ even lower than
the dimension of $[0,1]$, then we do get a positive anwer to our
question:

\begin{proposition}  \label{prop:Cantor}
Let $X$ be a totally disconnected space, and let $A$ be a
$C(X)$-algebra such that all fibres $A_x$, $x \in X$, of $A$ are
properly infinite. Then $A$ is properly infinite.
\end{proposition}

\begin{proof} Using Lemma~\ref{neighborhood} and the fact that $X$ is totally
  disconnected we can write $X$ as the disjoint union of clopen sets
  $F_1, F_2, \dots, F_n$ such that $A_{F_j}$ is properly infinite for
  all $j$. As 
$$A = A_{F_1} \oplus A_{F_2} \oplus \cdots \oplus A_{F_n},$$
the claim is proved.
\end{proof}

\vspace{.3cm}
\noindent{\sc Projet Algèbres d'opérateurs,
Institut de Mathématiques de Jussieu,
175, rue du Chevaleret, F-75013 PARIS, France}

\vspace{.2cm}
\noindent{\sl E-mail address:} {\tt  Etienne.Blanchard@math.jussieu.fr}\\
\noindent{\sl Internet home page:}
{\tt www.math.jussieu.fr/$\sim$blanchar}\\

\vspace{.5cm}

\noindent{\sc Department of Mathematics, University of Southern
  Denmark, Odense,
  Campusvej~55, 5230 Odense M, Denmark}

\vspace{.3cm}

\noindent{\sl E-mail address:} {\tt rohde@imada.sdu.dk} \\

\vspace{.5cm}

\noindent{\sc Department of Mathematics, University of Southern
  Denmark, Odense,
  Campusvej~55, 5230 Odense M, Denmark}

\vspace{.3cm}

\noindent{\sl E-mail address:} {\tt mikael@imada.sdu.dk}\\
\noindent{\sl Internet home page:}
{\tt www.imada.sdu.dk/$\sim$mikael/welcome} \\

\end{document}